\documentclass[11pt]{amsart}
\usepackage{amssymb}
\usepackage{url}

\newcommand{\base}{\operatorname{base}}
\newcommand{\roof}{\operatorname{roof}}
\newcommand{\down}{\operatorname{down}}
\newcommand{\up}{\operatorname{up}}

\newcommand{\A}{{\mathcal A}}
\newcommand{\Hecke}{{\mathcal H}}
\newcommand{\Z}{{\mathbb{Z}}}

\newcommand{\ulambda}{{\boldsymbol{\lambda}}}
\newcommand{\umu}{{\boldsymbol{\mu}}}
\newcommand{\res}{\operatorname{res}}
\newcommand{\unu}{{\boldsymbol{\nu}}}

\renewcommand{\leq}{\leqslant}
\renewcommand{\geq}{\geqslant}

\newtheorem{thm}{Theorem}[section]
\newtheorem{lem}[thm]{Lemma}
\newtheorem{cor}[thm]{Corollary}
\newtheorem{prop}[thm]{Proposition}
\newtheorem{conj}[thm]{Conjecture}

\theoremstyle{definition}
\newtheorem{exmp}[thm]{Example}
\newtheorem{defn}[thm]{Definition}

\theoremstyle{remark}
\newtheorem{rem}[thm]{Remark}

\begin{document}
\title[The Dipper-James-Murphy conjecture in
type ${\bf B_n}$]{Dipper-James-Murphy's conjecture for Hecke algebras of 
type ${\bf B_n}$}
\author{Susumu Ariki and Nicolas Jacon}

\address{S.A.: Research Institute for Mathematical Sciences, 
Kyoto University, Kyoto 606-8502, Japan}
\email{ariki@kurims.kyoto-u.ac.jp}

\address{N.J.:Universit\'e de Franche-Comt\'e,  UFR Sciences et
Techniques, 16 route de Gray, 25 030 Besan\c{c}on, France.}
\email{njacon@univ-fcomte.fr}

\date{March, 2007}

\begin{abstract}
We prove a conjecture by Dipper, James and Murphy that 
a bipartition is restricted if and only if it is Kleshchev. Hence 
the restricted bipartitions naturally label the crystal
graphs of level two irreducible integrable 
$\mathcal{U}_v({\widehat{\mathfrak{sl}}_e})$-modules 
and the simple modules of Hecke algebras of type $B_n$.
\end{abstract}

\maketitle

%


\centerline{\textit{Dedicated to Toshiaki Shoji and Ken-ichi Shinoda on their 
                      60$^{\text{th}}$ birthdays}}

\section{Introduction}

Let $F$ be a field, $q$ and $Q$ invertible elements of $F$. 
The Hecke algebra of type $B_n$ is the $F$-algebra 
defined by generators $T_0,\dots,T_{n-1}$ and relations 
\begin{gather*}
(T_0-Q)(T_0+1)=0,\quad (T_i-q)(T_i+1)=0\;(1\leq i<n)\\
(T_0T_1)^2=(T_1T_0)^2,\quad T_iT_{i+1}T_i=T_{i+1}T_iT_{i+1}\;(1\leq i<n-1)\\
T_iT_j=T_jT_i\;(j\ge i+2).
\end{gather*}
We denote it by $\Hecke_n(Q,q)$, or $\Hecke_n$ for short. 
The representation theory of $\Hecke_n$ in the semisimple case was studied by 
Hoefsmit, which had applications in determining generic degrees and 
Lusztig's $a$-values. Motivated by the modular representation theory of 
$U_n(q)$ and $Sp_{2n}(q)$ in the non-defining characteristic case, 
Dipper, James and Murphy began the study of 
the modular case more than a decade ago. The first task was 
to obtain classification of simple modules. For this, 
they constructed Specht modules which are indexed by the set of
bipartitions \cite{DJM}. 
The work shows in particular that Hecke algebras of 
type $B_n$ are cellular algebras in the sense of Graham and Lehrer. 
\footnote{This result has been recently generalized by Geck 
in \cite{Gcellular}.}
Then they conjectured that the simple modules
were labeled by $(Q,e)$-restricted bipartitions. 
Their philosophy to classify the simple $\Hecke_n$-modules resembles 
the highest weight theory in Lie theory:  
let $\mathfrak g$ be a semisimple Lie algebra. 
It has a commutative Lie subalgebra $\mathfrak h$, the Cartan subalgebra. 
One dimensional $\mathfrak h$-modules are called weights (by abuse of notion). 
When a $\mathfrak g$-module admits a simultaneous generalized eigenspace 
decomposition with respect to $\mathfrak h$, the decomposition is called 
the (generalized) weight space decomposition. 
Let $\Lambda$ be a weight. Suppose that a $\mathfrak g$-module $M$ has 
the property that 
\begin{itemize}
\item[(i)]
$\Lambda$ appears in the weight space decomposition of $M$,
\item[(ii)]
If $N$ is a proper $\mathfrak g$-submodule of $M$ then 
$\Lambda$ does not appear in the weight space decomposition of $N$.
\end{itemize}
Then the standard argument shows that $M$ has a unique nonzero 
irreducible quotient. In fact, Verma modules enjoy the property and 
their irreducible quotients give a complete set 
of simple objects in the BGG category. 
Now we turn to the Hecke algebra $\Hecke_n$. 
Define the Jucys-Murphy elements $t_1,\dots,t_n$ by 
$t_1=T_0$ and $t_{i+1}=q^{-1}T_it_iT_i$, for $1\leq i\leq n-1$. 
They generate a commutative subalgebra $\A_n$ of the Hecke algebra $\Hecke_n$, 
and $\A_n$ plays the role of the Cartan subalgebra: 
one dimensional $\A_n$-modules are called \emph{weights} and the generalized 
simultaneous eigenspace decomposition of an $\Hecke_n$-module is called 
the \emph{weight space decomposition}. 
Any weight is uniquely determined by the values at $t_1,\dots,t_n$ of 
the weight, 
and the sequence of these values in this order is called the 
\emph{residue sequence}.
Let $\ulambda=(\lambda^{(1)},\lambda^{(2)})$ be a bipartition (see \S 2.1) and 
let ${\bf t}$ be a standard bitableau of shape $\ulambda$ (see Def. \ref{tableau}). 
Then, ${\bf t}$ defines a weight whose values at $t_i$ are given by 
$c_iq^{b_i-a_i}$ where $a_i$ and $b_i$ are the row number and the 
column number of the node of ${\bf t}$ labelled by $i$ respectively, 
$c_i=-Q$ if the node is in $\lambda^{(1)}$ and 
$c_i=1$ if the node is in $\lambda^{(2)}$. 
By the theory of seminormal representations in the semisimple case 
and the modular reduction, a weight appears in 
some $\Hecke_n$-module if and only if it is obtained from 
a bitableau this way. 

Suppose that there is a weight obtained from a bitableau ${\bf t}$ of 
shape $\ulambda$ such that it does not appear in $S^\umu$ when 
$\umu\triangleleft\ulambda$. If such a bitableau exists, we say that 
$\ulambda$ is \emph{$(Q,e)$-restricted}. This is a clever generalization
of the notion of $e$-restrictedness.
Recall that a partition $\lambda=(\lambda_0,\lambda_1,\dots)$
is called \emph{$e$-restricted} if $\lambda_{i+1}-\lambda_{i}<e$, for all $i\ge0$.
Recall also that we have the similar Specht module theory for
Hecke algebras of type $A$.
Using Jucys-Murphy elements of the Hecke algebra of
type $A$, we can define weights as well. Then, a partition is
$e$-restricted if and only if
there is a weight obtained from a tableau of 
shape $\lambda$ such that it does not appear in $S^\mu$ when 
$\mu\triangleleft\lambda$. 

Recall from \cite{DJM} that 
$$
[S^\ulambda]=[D^\ulambda]+\sum_{\umu\triangleleft\ulambda}
d_{\ulambda\umu}[D^\umu],
$$
where the summation is over $\umu$ such that $D^\umu\ne0$, 
$d_{\ulambda\umu}$ are decomposition numbers, 
and $\sum_{\umu\triangleleft\ulambda}d_{\ulambda\umu}[D^\umu]$ is 
represented by the radical of the bilinear form on $S^\ulambda$. 
As $D^\umu$ is a surjective image of $S^\umu$, 
it implies that the weight does not appear in the 
radical, while it appears in $S^\ulambda$. Therefore, 
$D^\ulambda\ne0$ if $\ulambda$ is $(Q,e)$-restricted.
Unlike the case of the BGG category, we may have $D^\ulambda=0$ and 
it is important to know when it occurs. 
When $-Q$ is not a power of $q$, a bipartition 
$\ulambda=(\lambda^{(1)},\lambda^{(2)})$ 
is $(Q,e)$-restricted if and only if both $\lambda^{(1)}$ and 
$\lambda^{(2)}$ are $e$-restricted. 
Thus we know when a bipartition is 
$(Q,e)$-restricted. Further, \cite[Thm 4.18]{DJ} 
implies that $D^\ulambda\ne0$ if and only if $\ulambda$ is $(Q,e)$-restricted, 
that is, simple $\Hecke_n$-modules are labelled by $(Q,e)$-restricted bipartitions. 
Now we suppose that $-Q$ is a power of $q$. More precisely, we suppose that 
\begin{center}
\begin{itemize}
\item[(a)]
$q$ is a primitive $e^{th}$ root of unity with $e\geq2$,
\item[(b)]
$-Q=q^m$, for some $0\leq m<e$.
\end{itemize}
\end{center}
in the rest of the paper. We call $(Q,e)$-restricted bipartitions 
\emph{restricted} bipartitions. They conjectured in this case that 
$D^\ulambda\ne0$ if and only if $\ulambda$ is restricted, 
and it has been known as the Dipper-James-Murphy conjecture for Hecke 
algebras of type $B_n$. 

Later, connection with the theory of canonical bases in deformed Fock spaces 
in the sense of Hayashi and Misra-Miwa was discovered by 
Lascoux-Leclerc-Thibon \cite{LLT} and its proof in the framework of cyclotomic 
Hecke algebras \cite{A1} allowed the first author and Mathas \cite{A} 
\cite{AM} to label simple $\Hecke_n$-modules by the $n^{th}$ layer of
the crystal graph of the 
level two irreducible integrable $\mathfrak g(A^{(1)}_{e-1})$-module 
$L_v(\Lambda_0+\Lambda_m)$. In the theory, the crystal graph is 
realized as a subcrystal of the crystal of 
bipartitions, and the nodes of the crystal graph 
are called \emph{Kleshchev} bipartitions. More precise definition is
given in the next section and
$\ulambda=(\lambda^{(1)},\lambda^{(2)})$ is Kleshchev if and only if
$\lambda^{(2)}\otimes\lambda^{(1)}$
belongs to the subcrystal $B(\Lambda_0+\Lambda_m)$ of
$B(\Lambda_0)\otimes B(\Lambda_m)$, where the crystals $B(\Lambda_0)$ and
$B(\Lambda_m)$ are realized on the set of $e$-restricted partitions.
Now, $D^\ulambda\ne0$ if and only if $\ulambda$ is Kleshchev by \cite{A}. 
Hence, we obtained the classification of simple $\Hecke_n$-modules, 
or more precisely description of the set $\{\ulambda \mid D^\ulambda\ne0\}$, 
through a different approach and the Dipper-James-Murphy conjecture 
in the modern language is the statement that the Kleshchev bipartitions 
are precisely the restricted bipartitions.

The aim of this paper is to prove the Dipper-James-Murphy conjecture. 
Recall that Lascoux, Leclerc and Thibon considered Hecke algebras 
of type $A$ and they showed that if $\lambda$ is 
a $e$-restricted partition then we can find $a_1,\dots,a_p$ and 
$i_1,\dots,i_p$ such that we may write
$$
f^{(a_1)}_{i_1}...f^{(a_p)}_{i_p}{\emptyset}=
\lambda+\sum_{\nu\triangleright \lambda} c_{\nu,\lambda}(v) \nu
$$
in the deformed Fock space, where $c_{\nu,\lambda}(v)$ are Laurent 
polynomials. This follows from the ladder decomposition of a partition. 
Then LLT algorithm proves that Kleshchev partitions 
are precisely $e$-restricted partitions. 
The second author \cite{J} proved the similar formula for 
FLOTW multipartitions in the Jimbo-Misra-Miwa-Okado higher level 
Fock space using certain $a$-values instead of 
the dominance order. Recall that Geck and Rouquier gave
another method to label simple $\Hecke_n$-modules 
by bipartitions. The result shows that 
the parametrizing set of simple $\Hecke_n$-modules in the Geck-Rouquier theory,
which is called the \emph{canonical basic set}, 
is precisely the set of the FLOTW bipartitions. 
Our strategy to prove the conjecture is to give the analogous formula 
for Kleshchev bipartitions. To establish the formula, a non-recursive 
characterization of Kleshchev bipartitions given by the first author, 
Kreiman and Tsuchioka \cite{AKT} plays a key role. 

The paper is organized as follows. In the first section, 
we briefly recall the definition of Kleshchev bipartitions. 
We also recall the main result of \cite{AKT}. 
In the second section, we use this result to give an analogue for 
bipartitions of the ladder decomposition. 
Finally, the last section gives a proof for the conjecture.

{\bf Acknowledgments :} This paper was written when the second author 
visited RIMS in Kyoto. He would like to thank RIMS for 
the hospitality.

\section{Preliminaries} 

In this section, we recall the definition of Kleshchev bipartitions 
together with the main result of \cite{AKT} which gives a non-recursive 
characterization of these bipartitions. We fix $m$ as in the inroduction. 
Namely, the parameter $Q$ of the Hecke algebra is $Q=-q^m$ with $0\leq m<e$. 

\subsection{First definitions}\label{def} 
Recall that a partition $\lambda$ is a sequence of weakly decreasing 
nonnegative integers $(\lambda_0,\lambda_1,\cdots)$ such that 
$|\lambda|=\sum_{i\ge0}\lambda_i$ is finite. If $\lambda_i=0$ for
$i\ge r$ then we write $\lambda=(\lambda_0,\cdots,\lambda_{r-1})$.
A \emph{bipartition} $\ulambda$ is an ordered pair of partitions 
$(\lambda^{(1)},\lambda^{(2)})$. 
$|\ulambda|=|\lambda^{(1)}|+|\lambda^{(2)}|$ is called the \emph{rank} 
of $\ulambda$. The empty bipartition $(\emptyset,\emptyset)$ is the
only bipartition of rank zero.
The \emph{diagram} of $\ulambda$ is the set
$$\{ (a,b,c) \mid  1\leq c\leq 2, 0\leq b \leq \lambda^{(c)}_a-1\}\subseteq
\Z_{\ge0}^3.$$
We often identify a bipartition with its diagram. 
The \emph{nodes} of $\ulambda$ are the elements of the diagram. Let
$\gamma=(a,b,c)$ be a node of $\ulambda$. 
Then the \emph{residue} of $\gamma$ is defined by 
$$\res(\gamma)=\left\{ \begin{array}{ll}
b-a+m\;(\textrm{mod}\ e)\quad& \textrm{if }c=1, \\
b-a\;(\textrm{mod}\ e)\quad& \textrm{if }c=2. \\
\end{array}\right.$$
By assigning residues to the nodes of a bipartition, we view 
a bipartition as a colored diagram with colors in $\Z/e\Z$. 

\begin{exmp}\label{extableau}
Put $e=4$, $m=2$ and $\ulambda=((3,2),(4,2,1))$. Then the colored 
diagram associated with $\ulambda$ is as follows. 
$$
\left(
\begin{array}{|c|c|c|}
  \hline
  2 &3  &0     \\
  \hline
 1 & 2\\
\cline{1-2}
\end{array}\;,\;
\begin{array}{|c|c|c|c|}
  \hline
  0 &1  &2 & 3      \\
  \hline
 3 & 0\\
\cline{1-2}
2\\
 \cline{1-1}
\end{array}
\right)$$

\end{exmp}

If $\gamma$ is a node with residue $i$, we say that $\gamma$ is an \emph{$i$-node}.
Let $\ulambda$ and $\umu$ be two bipartitions such that 
$\umu=\ulambda\sqcup{\{\gamma\}}$. Then, we denote $\umu/\ulambda=\gamma$
and if $\res(\gamma)=i$, we say that $\gamma$ is a
\emph{removable $i$-node} of $\umu$. We also say that $\gamma$ is an
\emph{addable $i$-node} of $\ulambda$ by abuse of notion.
\footnote{An addable $i$-node of $\ulambda$ is not a node of $\ulambda$.}

Let $i\in\Z/e\Z$. We choose a total order on the set of removable and addable
$i$-nodes of a bipartition. Let $\gamma=(a,b,c)$ and 
$\gamma'=(a',b',c')$ be removable or addable $i$-nodes of a 
bipartition. We say that $\gamma$ is \emph{above} $\gamma'$ if either
$c=1$ and $c'=2$, or $c=c'$ and $a<a'$. \footnote{We now know that there are 
more than one Specht module theory, and different Specht module theories 
prefer different total orders on the set of $i$-nodes of a bipartition. 
Our choice of the total order is the one prefered by 
Dipper-James-Murphy's Specht module theory.}

Let $\mathcal F$ be the vector space over $\mathbb Q$ such that the basis
is given by the set of all bipartitions. We color the nodes of bipartitions as above.
We call it the (level two) \emph{Fock space}. 
We may equip it with $\hat{sl}_e$-module structure in which 
the action of the Chevalley generators is given by
$$
e_i\ulambda=\sum_{\unu:\res(\ulambda/\unu)=i}\unu,\quad
f_i\ulambda=\sum_{\unu:\res(\unu/\ulambda)=i}\unu.
$$
Using the total order on the set of removable and addable $i$-nodes given above, 
we deform the $\hat{sl}_e$-module structure to 
$U_v(\hat{sl}_e)$-module structure on the deformed Fock space 
$\mathcal F\otimes_{\mathbb Q}\mathbb Q(v)$, which is
the tensor product of two level one deformed Fock spaces.
We refer to \cite{Abook} for the details. 

\subsection{Kleshchev bipartitions}

Recall that the crystal basis of the deformed Fock space is given by the basis
vectors of the deformed Fock space. Hence it defines a crystal structure
on the set of bipartitions. We call it \emph{the crystal of bipartitions}.
As is explained in \cite{Abook}, the map
$(\lambda^{(1)},\lambda^{(2)})\mapsto \lambda^{(2)}\otimes\lambda^{(1)}$
identifies the crystal of bipartitions with the tensor product of
the crystal of partitions of highest weight $\Lambda_0$ and
that of highest weight $\Lambda_m$. As is already mentioned in the introduction,
Kleshchev bipartitions are those bipartitions 
which belongs to the same connected component as 
the empty bipartition in the crystal of bipartitions. Equivalently, 
Kleshchev bipartitions are those bipartitions which may be obtained from 
the empty bipartition by applying the Kashiwara operators successively. 
Rephrasing it in combinatorial terms, we have a recursive definition of 
Kleshchev bipartitions as follows. 

Let $\ulambda$ be a bipartition and let $\gamma$ be an $i$-node of
$\ulambda$, we say that $\gamma$ is a \emph{normal} $i$-node of $\ulambda$ if,
whenever $\eta$ is an addable $i$-node of $\ulambda$ below $\gamma$, 
there are more removable $i$-nodes between $\eta$
and $\gamma$ than addable $i$-nodes between $\eta$ and $\gamma$. If $\gamma$ is
the highest normal $i$-node of $\ulambda$, we say that $\gamma$ is a \emph{good} $i$-node of $\ulambda$. When $\gamma$ is a good $i$-node, we denote $\ulambda\setminus\{\gamma\}$ by $\tilde e_i\ulambda$.

\begin{defn}\label{Kleshchev}
A bipartition $\ulambda$ is \emph{Kleshchev} if either
$\ulambda=(\emptyset,\emptyset)$ or 
there exists $i\in\Z/e\Z$ and a good $i$-node $\gamma$ of $\ulambda$ such that
$\ulambda\setminus\{\gamma\}$ is Kleshchev.
\end{defn}
Note that the definition depends on $m$. The reader can prove easily using induction on $n=|\lambda^{(1)}|+|\lambda^{(2)}|$ that  if
$\ulambda=(\lambda^{(1)},\lambda^{(2)})$ is Kleshchev then 
both $\lambda^{(1)}$ and $\lambda^{(2)}$ are $e$-restricted. 
By general property of crystal bases, the following is clear.

\begin{lem}
Suppose that $\ulambda$ is a Kleshchev bipartition, $\gamma$ a good $i$-node of
$\ulambda$, for some $i$. Then $\tilde e_i\ulambda=\ulambda\setminus\{\gamma\}$ is Kleshchev.
\end{lem}

In \cite{AKT}, the first author, Kreiman and Tsuchioka have given a 
different characterization of Kleshchev bipartitions. 

Let $\lambda$ be a partition. Then the set of \emph{beta numbers of charge} 
$h$, where we only use $h=0$ or $h=m$ in the paper, is by definition
the set $J_h$ of decreasing integers
$$j_0>j_1>\cdots > j_k >\cdots $$
defined by $j_k=\lambda_k+h-k$, for $k\geq 0$. The charge $h$ also defines a coloring
of nodes: $\res(\gamma)=b-a+h\;(\textrm{mod}\ e)$ where
$a$ and $b$ are the row number and the 
column number of a node $\gamma$, respectively.

An addable $i$-node of $\lambda$ corresponds to $x\in J_h$ such that
$x+e\Z=i$ and $x+1\not\in J_h$. We call $x$ an addable $i$-node of $J_h$.
Similarly, a removable $i$-node of $\lambda$ corresponds to $x\in J_h$ such that
$x+e\Z=i+1$ and $x-1\not\in J_h$. We call $x$ a removable $i$-node of $J_h$.

We define the abacus display of $J_h$ in the usual way.
The $i^{th}$ runner of the abacus is
$\{x\in\Z \mid x+e\Z=i\}$, for $i\in\Z/e\Z$.

\begin{defn}
Let $\lambda$ be an $e$-restricted partition, $J_m(\lambda)$ 
the corresponding set of beta numbers of charge $m$. We write $J_m$ for
$J_m(\lambda)$ and define 
$$U(J_m)=\{x\in J_m \mid x-e\notin J_m \}.$$
If $\lambda$ is an $e$-core then we define $\up_m(\lambda)=\lambda$. 
Otherwise let $p=\max U(J_m)$ and define 
$$V(J_m)=\{x>p \mid x\neq p\; (\textrm{mod}\ e),x-e\in J_m ,x\notin J_m \}.$$
Note that $V(J_m)$ is nonempty since $\lambda$ is $e$-restricted. Let $q=\min V(J_m)$. Then we define 
$$\up(J_m)=(J_m \setminus \{p\})\sqcup \{q\}$$ and we denote the 
corresponding partition by $\up_m(\lambda)$. 
\end{defn}

In \cite{AKT}, it is shown that $\up_m(\lambda)$ is again $e$-restricted and
we reach an $e$-core after applying $\up_m$ finitely many times.

\begin{defn}
Let $\lambda$ be an $e$-restricted partition. 
Apply $\up_m$ repeatedly until we reach an $e$-core. 
We denote the resulting $e$-core by $\roof_m{(\lambda)}$.
\end{defn}

\begin{defn}
Let $\lambda$ be an $e$-restricted partition, $J_0(\lambda)$ 
the corresponding set of beta numbers of charge $0$. We write $J_0$ for
$J_0(\lambda)$ and define 
$$U(J_0)=\{x\in J_0 \mid x-e\notin J_0\}.$$
If $\lambda$ is an $e$-core then we define $\down_0(\lambda)=\lambda$. 
Otherwise let $p'=\min U(J_0)$ and define 
$$W(J_0)=\{x>p'-e \mid x\in J_0, x+e\notin J_0\}\cup\{p'\}.$$
It is clear that $W(J_0)$ is nonempty. Let $q'=\min W(J_0)$. Then we define 
$$\down (J_0)=(J_0\setminus \{q'\})\sqcup \{p'-e\}$$ and we denote 
the corresponding partition by $\down_0(\lambda)$. 
\end{defn}

In \cite{AKT}, it is shown that $\down_0(\lambda)$ is again $e$-restricted and
we reach an $e$-core after applying $\down_0$ finitely many times.

\begin{defn}
Let $\lambda$ be an $e$-restricted partition. 
Apply $\down_0$ repeatedly until we reach an $e$-core. 
We denote the resulting $e$-core by $\base_0{(\lambda)}$.
\end{defn}

Finally, let $\lambda$ be an $e$-restricted partition, 
$J^{\max}_0$ the set of beta numbers of charge $0$ for 
$\base_0(\lambda)$. Define 
$M_i(\lambda)$, for $i\in\Z/e\Z$, by
$$M_i(\lambda)=\max\{x\in J^{\max}_0 \mid x+e\Z=i\}.$$
We write $M_i(\lambda)$ in decreasing order
$$M_{i_1}(\lambda) >M_{i_2}(\lambda) >\cdots > M_{i_e}(\lambda).$$
Then $J^{\max}_0\cup \{M_{i_k}(\lambda)+e\}_{1\leq k\leq m}$ is 
the set of beta numbers of charge $m$, for some partition. 
We denote the partition by $\tau_m(\base_0(\lambda))$. 

Now, the characterization of Kleshchev bipartitions is as follows. 

\begin{thm}[\cite{AKT}] Let $\ulambda=(\lambda^{(1)},\lambda^{(2)})$ be a 
bipartition such that both $\lambda^{(1)}$ and $\lambda^{(2)}$ are 
$e$-restricted. Then $\ulambda$ is Kleshchev if and only if
$$\roof_m(\lambda^{(1)})\subseteq \tau_m(\base_0(\lambda^{(2)})).$$
\end{thm}

\subsection{The Dipper-James-Murphy conjecture}
We recall the dominance order for bipartition. Let 
$\ulambda=(\lambda^{(1)},\lambda^{(2)})$ and $\umu=(\mu^{(1)},\mu^{(2)})$ 
be bipartitions. In this paper, we write 
$\umu\trianglelefteq\ulambda$ if 
$$
\sum_{k=1}^j \lambda^{(1)}_k\geq \sum_{k=1}^j \mu^{(1)}_k \;\text{and}\;
|\lambda^{(1)}|+\sum_{k=1}^j \lambda^{(2)}_k
\geq |\mu^{(1)}|+\sum_{k=1}^j \mu^{(2)}_k,
$$
for all $j\geq0$. 

\begin{defn}\label{tableau}
Let $\ulambda$ be a bipartition of rank $n$. 
A \emph{standard bitableau of shape $\ulambda$} is a 
sequence of bipartitions 
$$
{\boldsymbol\emptyset}=\ulambda[0]\subseteq\ulambda[1]\subseteq\cdots
\subseteq\ulambda[n]=\ulambda
$$
such that the rank of $\ulambda[k]$ is $k$, for $0\leq k\leq n$. 
Let ${\bf t}$ be a standard bitableau of shape $\ulambda$. Then
the \emph{residue sequence} of ${\bf t}$ is the sequence
$$
(\res(\gamma[1]),\dots,\res(\gamma[n]))\in (\Z/e\Z)^n
$$
where $\gamma[k]=\ulambda[k]/\ulambda[k-1]$, for $1\leq k\leq n$. 
\end{defn}

A standard bitableau may be viewed as filling of the nodes of $\ulambda$ 
with numbers $1,\dots,n$: we write $k$ in the node $\gamma[k]$, for 
$1\leq k\leq n$. 

\begin{defn}\label{rest} A bipartition $\ulambda$ is 
\emph{$(-q^m,e)$-restricted}, or \emph{restricted} for short, if 
there exists a standard bitableau ${\bf t}$ of shape $\ulambda$ 
such that the residue sequence of any standard bitableau of shape 
$\unu\lhd\ulambda$ does not coincide with the residue sequence of ${\bf t}$. 
\end{defn}

\begin{conj}[{\cite[Conj. 8.13]{DJM}}]
A bipartition $\ulambda$ is Kleshchev if and only if it is restricted. 
\end{conj}

\section{Properties of Kleshchev bipartitions}
The aim of this section is to prove some combinatorial results 
concerning Kleshchev bipartitions.

\subsection{Admissible sequence} 

\begin{defn}
Let $i\in\Z/e\Z$. We say that 
a sequence of removable $i$-nodes $R_1,....,R_s$ (where $s\ge 1$)
of a bipartition $\ulambda$ is an
\emph{admissible sequence of $i$-nodes for $\ulambda$} if
\begin{itemize}
\item $R_1,....,R_s$ are the lowest $s$
 removable $i$-nodes of $\ulambda$ and every addable $i$-nodes is
 above all of these nodes, and
\item if there is a removable $i$-node $R$ above $R_1,....,R_s$, 
there must exist an addable $i$-node below $R$.
\end{itemize}
\end{defn}

The following lemma shows the existence of an admissible sequence 
of $i$-nodes, for some $i$, for a Kleshchev bipartition: choose $i$ as in the lemma and 
read addable and removable $i$-nodes in the total order of nodes. 
Suppose that $\ulambda$ has at least one addable $i$-node and
let $\eta$ be the lowest addable $i$-node. Then removable $i$-nodes 
below $\eta$ form an admissible sequence of $i$-nodes.
If $\ulambda$ does not have an addable $i$-node,
all removable $i$-nodes of $\ulambda$ form an admissible sequence of $i$-nodes.

\begin{lem}\label{Klesh} Let $\ulambda=(\lambda^{(1)},\lambda^{(2)})$ 
be a nonempty Kleshchev bipartition. 
Then there exists $i\in\Z/e\Z$ and a removable $i$-node $\gamma$ 
such that if $\eta$ is an addable $i$-node of $\ulambda$ 
then $\eta$ is above $\gamma$.
\end{lem}
\begin{proof}
Recall that both $\lambda^{(1)}$ and $\lambda^{(2)}$ are $e$-restricted. 
There are two cases to consider.
\begin{itemize}
\item 
Assume that $\lambda^{(2)}$ is not the empty partition. 
Let $\lambda^{(2)}_j$ be the last row. Define 
$\gamma=(j,\lambda^{(2)}_j-1,2)$ and $i=\res(\gamma)$. 
Since $\lambda^{(2)}$ is $e$-restricted, 
the residue of the addable node $(j+1,0,2)$ is not $i$. 
Hence all the addable $i$-node of $\ulambda$ are above $\gamma$.
\item 
Assume that $\lambda^{(2)}$ is the empty partition. Let $\lambda^{(1)}_j$ be 
the last row. Define $\gamma=(j,\lambda^{(1)}_j-1,1)$ and $i=\res(\gamma)$. 
Since $\lambda^{(1)}$ is $e$-restricted, the residue of the addable node 
$(j+1,0,1)$ is not $i$. We show that the residue of the addable 
node $(0,0,2)$ is not $i$. Suppose to the contrary that the residue is $i$. 
As $\ulambda$ is Kleshchev, we may delete good nodes successively to obtain 
the empty bipartition. Hence $\gamma$ must be deleted 
at some point in the process. However, it can never be a good node since 
we always have an addable $i$-node $(0,0,2)$ just below it and there is no 
removable $i$-node between them, so we have a contradiction. 
\end{itemize}
Hence the claim follows.
\end{proof}

\begin{lem} Let $\ulambda=(\lambda^{(1)},\lambda^{(2)})$ be a
nonempty Kleshchev bipartition and 
$R_1,....,R_s$ an admissible sequence of $i$-nodes
for $\ulambda$.
Define $\umu=(\mu^{(1)},\mu^{(2)})$ by $\ulambda=\umu\sqcup\{R_1,....,R_s\}$.
Then $\umu$ is Kleshchev.
\end{lem}
\begin{proof}
Recall that $\lambda^{(1)}$ and $\lambda^{(2)}$ are both $e$-restricted.
We claim that $\mu^{(1)}$ and $\mu^{(2)}$ are both $e$-restricted. 
We only prove that $\mu^{(1)}$ is $e$-restricted as the proof for 
$\mu^{(2)}$ is the same. Suppose that $\mu^{(1)}$ is not $e$-restricted.
Since $\lambda^{(1)}$ is $e$-restricted, 
it occurs only when there exists 
$j$ such that $\lambda_j^{(1)}=\lambda^{(1)}_{j+1}+e-1$,
$\mu_j^{(1)}=\lambda_j^{(1)}$, $\mu_{j+1}^{(1)}=\lambda_{j+1}^{(1)}-1$ 
and $\res(j+1,\lambda_{j+1}^{(1)}-1,1)=i$. But then 
$\res(j,\lambda_j^{(1)}-1,1)=i$, which implies $\mu_j^{(1)}=\lambda_j^{(1)}-1$ 
by definition of $\umu$. \\
\\
{\bf First case.} First we consider the case when either
$\lambda^{(2)}=\emptyset$ or $\lambda^{(2)}\ne\emptyset$ and
$\lambda^{(2)}$ has no addable $i$-node. 
If $\lambda^{(1)}$ has no addable 
$i$-node then the admissible sequence $R_1,....,R_s$ is given by all the
removable $i$-nodes of $\ulambda$, and thus
all the normal $i$-nodes of $\ulambda$.
Hence $\umu=\tilde e_i^s\ulambda$, which implies that $\umu$ is Kleshchev.
Therefore, we may and do assume that $\lambda^{(1)}$ has at least one addable $i$-node.

Define $t\ge0$ by
$$
t=\min\{k \mid \roof_m(\lambda^{(1)})=\up_m^k(\lambda^{(1)})\}.
$$
We prove that $\umu$ is Kleshchev by induction on $t$. Note that 
if $\ulambda$ is Kleshchev then so is 
$(\up_m(\lambda^{(1)}),\lambda^{(2)})$ since 
$$
\roof_m(\up_m(\lambda^{(1)}))=\roof_m(\lambda^{(1)})\subseteq
\tau_m(\base_0(\lambda^{(2)})).
$$

Suppose that $t=0$. Then $\lambda^{(1)}$ is an $e$-core. 
As $\lambda^{(1)}$ has an addable $i$-node,
$\lambda^{(1)}$ has no removable $i$-node. Thus
all the removable $i$-nodes of
$\ulambda$ are nodes of $\lambda^{(2)}$.
As $\lambda^{(2)}$ has no addable
$i$-node, the admissible sequence $R_1,....,R_s$ is
given by all the normal $i$-nodes of $\ulambda$. Hence,
$\umu=\tilde e_i^s\ulambda$ and $\umu$ is Kleshchev. 

Suppose that $t>0$ and that the lemma holds for 
$(\up_m(\lambda^{(1)}),\lambda^{(2)})$. Recall that we have assumed that
$\lambda^{(1)}$ has an addable $i$-node. Let $r$ be the 
minimal addable $i$-node of $J_m:=J_m(\lambda^{(1)})$.
The corresponding addable $i$-node of $\lambda^{(1)}$, say $\gamma$,
is the lowest addable $i$-node of $\ulambda$ and the admissible sequence
$R_1,....,R_s$ is given by all the removable $i$-nodes of $\ulambda$ that is
below $\gamma$. If there is no removable $i$-node 
greater than $r$ then all the removable $i$-nodes of $\ulambda$ 
are normal and by deleting them, we obtain that $\umu$ is 
Kleshchev. Hence, we assume that there is a removable $i$-node 
greater than $r$.
As $r+1\not\in J_m$, this implies that there is $x\in U(J_m)$
on the $(i+1)^{th}$ runner such that $x>r+1$. 
Let $p=\max U(J_m)$. Then $x\in U(J_m)$ implies that 
$p\ge x>r+1$.
As $p$ moves to $q>p$,
it implies that $R_1,....,R_s$ is an admissible sequence of $i$-nodes for
$(\up_m(\lambda^{(1)}),\lambda^{(2)})$ and that
$$ (\up_m(\lambda^{(1)}),\lambda^{(2)})=
(\up_m(\mu^{(1)}),\mu^{(2)})\sqcup \{R_1,....,R_s\}. $$
Now, $(\up_m(\mu^{(1)}),\mu^{(2)})$  is Kleshchev by the induction 
hypothesis. Hence,
$$
\roof_m(\mu^{(1)})=\roof_m(\up_m(\mu^{(1)}))\subseteq\tau_m(\base_0(\mu^{(2)}))
$$
and $\umu$ is Kleshchev as desired.\\ 
\\
{\bf Second case :} Now, we consider the case when $\lambda^{(2)}\ne\emptyset$
and $\lambda^{(2)}$ has at least one addable $i$-node. 
Note that it forces $\lambda^{(1)}=\mu^{(1)}$ and
$R_1,....,R_s$ are nodes of $\lambda^{(2)}$. Define $t'\ge0$ by
$$
t'=\min\{k \mid \base_0(\lambda^{(2)})=\down_0^k(\lambda^{(2)})\}.
$$
We prove that $\umu$ is Kleshchev by induction on $t'$. Note that 
if $\ulambda$ is Kleshchev then so is 
$(\lambda^{(1)},\down_0(\lambda^{(2)}))$ since
$$
\roof_m(\lambda^{(1)})\subseteq\tau_m(\base_0(\lambda^{(2)}))
=\tau_m(\base_0(\down_0(\lambda^{(2)}))).
$$

If $t'=0$ then $\lambda^{(2)}$ 
is an $e$-core and it has removable $i$-nodes $R_1,....,R_s$. 
Hence, $\lambda^{(2)}$ has no addable $i$-node and we are reduced to 
the previous case. Thus we suppose that $t'>0$ and that 
the lemma holds for $(\lambda^{(1)},\down(\lambda^{(2)}))$.  
Let $J=J_0 (\lambda^{(2)})$ and
$$
r=\min\{x\in J \mid x+e\Z=i+1,\;x-1\notin J \}-1.
$$
Note that $r$ is on the $i^{th}$ runner. Then there exists $N\ge1$ such that 
$$r,r+e,\dots,r+(N-1)e\not\in J \quad\text{and}\quad r+Ne\in J .$$

Let $p'=\min U(J)$. Then $p'\le r+Ne$. 
Suppose that $p'$ is not on the $i^{th}$ runner or   
the $(i+1)^{th}$ runner. If a node which is not on one of these two runners 
moves to $p'-e$ by the down operation, the admissible sequence
$R_1,....,R_s$ is an admissible sequence of
$i$-nodes for $(\lambda^{(1)},\down_0(\lambda^{(2)}))$ and
$$(\lambda^{(1)},\down_0(\lambda^{(2)}))=
(\mu^{(1)},\down_0 (\mu^{(2)}))\sqcup \{R_1,....,R_s\}.$$
Thus, by the induction hypothesis, 
$(\mu^{(1)},\down_0(\mu^{(2)}))$ is Kleshchev. Hence,
$$
\roof_m(\mu^{(1)})\subseteq\tau_m(\base_0(\down_0(\mu^{(2)})))
=\tau_m(\base_0(\mu^{(2)}))
$$
implies that $\umu$ is Kleshchev.

If a node in one of the two runners 
moves to $p'-e$, then there exists $0\le k\le N-1$ such that 
$r+ke+1\in J$, $r+(k+1)e+1\not\in J $ and $r+ke+1$ moves 
to $p'-e$. Suppose that $k<N-1$. 
Then, $r+ke\in J_0(\mu^{(2)})$ and $r+(k+1)e\not\in J_0(\mu^{(2)})$.
Hence, $r+ke$ moves to $p'-e$ to obtain $\down_0(\mu^{(2)})$.
As $r+ke+1$ corresponds to one of $R_1,....,R_s$, say $R_k$,
$R_1,\dots,\hat{R_k},\dots,R_s$ is an admissible sequence of $i$-nodes
for $(\lambda^{(1)},\down_0(\lambda^{(2)}))$ and
$$(\lambda^{(1)},\down_0(\lambda^{(2)}))=
(\mu^{(1)},\down_0(\mu^{(2)}))\sqcup \{R_1,\dots,\hat{R_k},\dots,R_s\}.$$
Hence, $(\mu^{(1)},\down_0(\mu^{(2)}))$ is Kleshchev
by the induction hypothesis, and $\umu$ is Kleshchev as before.
Next suppose that $k=N-1$. As $r+(N-1)e+1$ moves to $p'-e$, we have
$$r+(N-1)e+1<p'<r+Ne,$$
$r+Ne+1\notin J_0$ and $r+Ne$ is an addable $i$-node. Let 
$K$ be the set of beta numbers of charge $0$ of 
$\mu^{(2)}$. For $x\in \mathbb{Z}$, we denote $J_{\le x}:=J\cap \mathbb{Z}_{\leq x}$
and $K_{\le x}:=K\cap \mathbb{Z}_{\leq x}$.
We claim that 
$$
\base(J_{\le r+Ne})=\base(K_{\le r+Ne}).
$$
Let $p'=y_0<y_1<\cdots<y_l<r+Ne$ be the nodes in $J$ which are 
greater than or equal to $p'$ and smaller than $r+Ne$. 
We show that 
$$\base(J_{\le y_j})=s_i\base(K_{\le y_j})\supseteq\base(K_{\le y_j}),$$
for $0\le j\le l$, where $s_i$ means swap of the $i^{th}$ and 
$(i+1)^{th}$ runners.  
$J_{\le p'}$ and $K_{\le p'}$ are $s_i$-cores in the sense of \cite{AKT}, 
and direct computation shows the formula for $j=0$. Now we use 
$\base(J_{\le y_{j+1}})=\base(\{y_{j+1}\}\cup\base(J_{\le y_j}))$
and $\base(K_{\le y_{j+1}})=\base(\{y_{j+1}\}\cup\base(K_{\le y_j}))$ 
\footnote{See \cite[Prop 7.8]{AKT}.} and 
continue the similar computation and comparison of 
$\base(J_{\le y_j})$ and $\base(K_{\le y_j})$. At the end of the 
inductive step, we get
$$
\base(J_{<r+Ne})=s_i\base(K_{<r+Ne})\supseteq\base(K_{<r+Ne}).
$$
Now, one more direct computation shows
$$
\base(\{r+Ne\}\cup\base(J_{<r+Ne}))=\base(\{r+Ne\}\cup\base(K_{<r+Ne})),
$$
and we have the claim. Therefore, $\base_0(\lambda^{(2)})=\base_0(\mu^{(2)})$ 
and we have 
$$
\roof_m(\lambda^{(1)})\subseteq\tau_m(\base_0(\lambda^{(2)}))
=\tau_m(\base_0(\mu^{(2)})).
$$
Recalling that $\lambda^{(1)}=\mu^{(1)}$, we have that $\umu$ is Kleshchev.

Suppose that $p'$ is on one of the two runners. As $p'\le r+Ne$, we 
have either $p'=r+Ne$ or $p'=r+ke+1$, for some $0\le k\le N-1$. If 
the latter occurs, $\down_0(\lambda^{(2)})$ is obtained by 
moving a node outside the two runners to $p'-e$ or moving 
$p'$ to $p'-e$, and $\down_0(\mu^{(2)})$ is 
obtained from $\mu^{(2)}$ by moving the same node outside the 
two runners to $p'-1-e$ or moving $p'-1$ to $p'-1-e$, respectively.
Thus, $\down_0(\mu^{(2)})$ is obtained from $\down_0(\lambda^{(2)})$
by removing the nodes of an admissible sequence of
$i$-nodes for $(\lambda^{(1)},\down_0(\lambda^{(2)}))$.
Hence $(\mu^{(1)},\down_0(\mu^{(2)}))$ is Kleshchev
by the induction hypothesis, and it follows that $\umu$ is Kleshchev.
If $p'=r+Ne$ and $r+Ne+1\in J_0$ then the same is true, and 
if $p'=r+Ne$ and $r+Ne+1\not\in J_0$ then $\mu^{(2)}=\down_0^N(\lambda^{(2)})$ 
and we have 
$$
\roof_m(\lambda^{(1)})\subseteq\tau_m(\base_0(\lambda^{(2)}))
=\tau_m(\base_0(\mu^{(2)})).
$$
Hence, $\umu=(\lambda^{(1)},\mu^{(2)})$ is Kleshchev. 
\end{proof}

\subsection{Optimal sequence of a Kleshchev bipartition} 
Let $\ulambda=(\lambda^{(1)},\lambda^{(2)})$ be a Kleshchev bipartition. 
By the previous lemma, we may define by induction a sequence of Kleshchev 
bipartitions 
$$
\ulambda=:\ulambda[1],\;\ulambda[2],...,\;\ulambda[p],\;\ulambda[p+1]=
{\boldsymbol\emptyset}
$$
and a sequence of residues 
$$
\underbrace{i_1,...,i_1}_{a_1\textrm{ times}},...,
\underbrace{i_p,...,i_p}_{a_p\textrm{ times}}
$$
such that, for $1\le j\le p$, $\ulambda[j]=\ulambda[j+1]\sqcup\{ R^j_1,....,R_{a_j}^j\}$
and $R^j_1,....,R^j_{a_j}$ is an admissible sequence of $i_j$-nodes
for $\ulambda[j]$. 

We call
$\underbrace{i_1,...,i_1}_{a_1\textrm{ times}},...,
\underbrace{i_p,...,i_p}_{a_p\textrm{ times}}$ an \emph{optimal sequence} 
of $\ulambda$.

\begin{exmp} Keeping example \ref{extableau}, it is easy to see that 
$((3,2),(4,2,1))$ is a Kleshchev bipartition and an optimal 
sequence is given by
$$2,2,0,3,3,2,1,1,0,0,3,2.$$
\end{exmp}

\section{Proof of the conjecture}
\subsection{The result}
We are now ready to prove the conjecture. As is explained in the introduction, 
it is enough to prove that Kleshchev bipartitions are restricted bipartitions. 
To do this, we introduce certain reverse lexicographic order on the set of
bipartitions. 

\begin{defn}
We write $\ulambda \prec \unu$ if either
\begin{itemize}
\item
there exists $j\ge 0$ such that $\lambda^{(2)}_k=\nu^{(2)}_k$, for $k>j$, and $\lambda^{(2)}_j< \nu^{(2)}_j$, or
\item 
there exists $j\ge 0$ such that $\lambda^{(2)}=\nu^{(2)}$, 
$\lambda^{(1)}_k=\nu^{(1)}_k$, for $k>j$, and $\lambda^{(1)}_j< \nu^{(1)}_j$.
\end{itemize}
\end{defn}
It is clear that if $\unu \lhd \ulambda$ then $\ulambda\prec \unu$. 
Recall that the deformed Fock space is given a specific 
$U_v(\hat{sl}_e)$-module structure which is suitable for the 
Dipper-James-Murphy's Specht module theory. 

\begin{prop}\label{DJM1}
  Let $\ulambda=(\lambda^{(1)},\lambda^{(2)})$ be Kleshchev and let 
$$\underbrace{i_1,...,i_1}_{a_1\textrm{ times}},...,
\underbrace{i_p,...,i_p}_{a_p\textrm{ times}}$$
be an optimal sequence of $\ulambda$. Then we have
$$
f^{(a_1)}_{i_1}...f^{(a_p)}_{i_p}{\boldsymbol\emptyset}=
\ulambda+\sum_{\unu\prec \ulambda} c_{\unu,\ulambda}(v) \unu,
$$
for some Laurent polynomials $c_{\unu,\ulambda}(v)\in\Z_{\ge0}[v,v^{-1}]$, in the 
deformed Fock space. 
\end{prop}
\begin{proof} First note that the coefficient of $\ulambda$ is one because
each admissible sequence of $i_j$-nodes is a sequence of normal $i_j$-nodes.

Now the proposition is proved by induction on $p$. Let $R^1_1,....,R^1_{a_1}$ be 
the admissible sequence of $i_1$-nodes for $\ulambda$ and let 
$\ulambda'$ be the Kleshchev bipartition such that 
$$\ulambda=\ulambda'\sqcup\{ R^1_1,....,R_{a_1}^1\}.$$
By the induction hypothesis, we have 
$$f^{(a_2)}_{i_2}...f^{(a_{p})}_{i_{p}}{\bf \emptyset}=\ulambda'+\sum_{\unu'\prec \ulambda'} c_{\unu,\ulambda'}(v) \unu'.$$
Let $\unu\neq \ulambda$ be a bipartition such that it appears in 
$f^{(a_1)}_{i_1}...f^{(a_p)}_{i_p}{\bf \emptyset}$ with nonzero coefficient. 
Then there exist removable $i_1$-nodes 
${R'}_1^1 $,...,${R'}^1_{a_1}$ of $\unu$ and
a bipartition $\unu'\preceq \ulambda'$ such that 
$$\unu=\unu'\sqcup\{ {R'}^1_1,....,{R'}_{a_1}^1\}.$$
As $R^1_1,....,R_{a_1}^1$ are the lowest $a_1$ $i_1$-nodes of $\ulambda$,
$\unu'=\ulambda'$ implies $\unu\preceq \ulambda$. Hence 
we may assume $\unu'{\prec} \ulambda'$. 
Suppose that we have $\ulambda\prec \unu$. 
If ${\nu'}^{(2)}\ne{\lambda'}^{(2)}$ then choose $t$ such that
${\nu'}^{(2)}_t<{\lambda'}^{(2)}_t$ and
${\nu'}^{(2)}_j={\lambda'}^{(2)}_j$, for $j>t$. Then we can show
\begin{itemize}
\item[(i)]
$\nu^{(2)}_j=\lambda^{(2)}_j$, for $j>t$, 
\item[(ii)]
$\nu^{(2)}_{t+1}<\nu^{(2)}_t={\nu'}^{(2)}_t+1={\lambda'}^{(2)}_t=\lambda^{(2)}_t$,
\item[(iii)]
at least one of the nodes $R^1_1,....,R_{a_1}^1$ 
is above $(t,\lambda^{(2)}_t-1,2)$. 
\end{itemize}
The condition (ii) implies that 
$\res(t,\lambda^{(2)}_t-1,2)=\res(t,{\nu'}^{(2)}_t,2)=i_1$.
Thus $(t,\lambda^{(2)}_t-1,2)$ is an $i_1$-node
of $\lambda^{(2)}$. Then (iii) implies that 
it is not a removable node of $\lambda^{(2)}$
(otherwise it has to be removed to 
obtain $\ulambda'$). This implies that 
${\lambda}^{(2)}_{t}={\lambda}^{(2)}_{t+1}$. 
Thus $\nu^{(2)}_{t+1}<\lambda^{(2)}_t=\lambda^{(2)}_{t+1}$
and (i) is contradicted.

If ${\nu'}^{(2)}={\lambda'}^{(2)}$ then choose $t$ such that
${\nu'}^{(1)}_t<{\lambda'}^{(1)}_t$ and
${\nu'}^{(1)}_j={\lambda'}^{(1)}_j$, for $j>t$. Then we argue as above.
\end{proof}

\begin{cor}
The Dipper-James-Murphy conjecture is true.
\end{cor}
\begin{proof}
Observe that $\unu$ appears in $f_{i_n}\cdots f_{i_1}{\boldsymbol\emptyset}$ 
if and only if there exists a standard bitableau of shape $\unu$ 
such that its residue sequence is $(i_1,\dots,i_n)$. 
Let $\ulambda$ be Kleshchev. Then Proposition \ref{DJM1} shows that 
there is a standard bitableau ${\bf t}$ such that if the residue sequence 
of ${\bf t}$ appears as the residue sequence of a standard bitableau 
of shape $\unu$ then $\unu\preceq \ulambda$. 
Suppose that the residue sequence of ${\bf t}$ is the residue sequence 
of a standard bitableau of shape $\unu\triangleleft\ulambda$. 
As $\unu\triangleleft\ulambda$ implies $\ulambda\prec \unu$, we cannot 
have $\unu\preceq\ulambda$, a contradiction. 
Hence $\ulambda$ is restricted. 
\end{proof}

\subsection{Remark} 
We conclude the paper with a remark. 

\begin{rem}
There is a systematic way to produce realizations of
the crystal $B(\Lambda_0+\Lambda_m)$ on 
a set of bipartitions, for each choice of $\log_q(-Q)$.
The bipartitions are called Uglov bipartitions. 
Recent results of Geck \cite{Gcellular} and Geck and the second author 
\cite{GJ} show that Uglov bipartitions naturally label simple 
$\Hecke_n$-modules, and Rouquier's theory of 
the BGG category of rational Cherednik algebras as 
quasihereditary covers of Hecke algebras naturally explains the 
existence of various Specht module theories which depends on $\log_q(-Q)$. 

We conjecture that Uglov bipartitions satisfy an analogue of 
Proposition \ref{DJM1} except that the dominance order is 
replaced by an appropriate 
$a$-value in the sense of \cite[Prop 2.1]{GJ}. As is mentioned in the 
introduction, it is known that 
this conjecture is true in the case where Uglov bipartitions are FLOTW 
bipartitions \cite[Prop. 4.6]{J}.
\end{rem}



\begin{thebibliography}{131}

\bibitem{A1}
{\sc S.~Ariki}, On the decomposition numbers of the Hecke algebra of 
$G(m,1,n)$, J.~Math.~Kyoto Univ., \textbf{36} (1996), 789--808.

\bibitem{Abook}
{\sc S.~Ariki}, Representations of quantum algebras and combinatorics of Young
tableaux, Univ. Lecture Series, \textbf{26} (2002), AMS.

\bibitem{A}
{\sc S.~Ariki}, On the classification of simple modules for cyclotomic 
Hecke algebras of type $G(m,1,n)$ and Kleshchev multipartitions, 
Osaka J. Math.  38  (2001),  no. 4, 827--837. 

\bibitem{AKT}
{\sc S.~Ariki, V.~Kreiman and S.Tsuchioka}, On the tensor product of two basic
representations of $U_v(\widehat{sl}_e)$, to appear in Advances in Math. (available at {\tt
http://arXiv.org/math.RT/0606044}.)

\bibitem{AM}
{\sc S.~Ariki and A~Mathas}, The number of simple modules of the Hecke algebras 
of type G(r,1,n), Math.~Zeitschrift, 233 (2000), 601--623.

\bibitem{DJ}
{\sc R.~Dipper and G. D.~James},  Representations of Hecke algebras of type $B\sb n$, 
J.~Algebra {\bf 146} (1992), 454--481. 

\bibitem{DJM}
{\sc R.~Dipper, G. D.~James and G. E.~Murphy}, Hecke algebras of type $B_n$ at
roots of unity, Proc. London Math. Soc. {\bf 70} (1995), 505--528.

\bibitem{Gcellular}
{\sc M.~Geck}, Hecke algebras of finite type are cellular, Invent. Math. {\bf 169}, Vol. 3, (2007), 501-517.

\bibitem{GJ}
{\sc M.~Geck and N. ~Jacon}, Canonical basic sets in type $B_n$, 
J. Algebra 306 (2006), 104-127.

\bibitem{J}
{\sc N. ~Jacon}, 
On the parametrization of the simple modules for Ariki-Koike algebras at 
roots of unity, J. Math. Kyoto Univ. 44 (2004), no. 4, 729-767.

\bibitem{LLT}
{\sc A.~Lasoux, B.~Leclerc and J-Y.~Thibon}, Hecke algebras at roots
of unity and cystal bases of quantum affine algebras, Comm. Math.
Phys. {\bf 181} (1996), 205--263.

\bibitem{M}
{\sc A.~Mathas}, The representation theory of the Ariki-Koike and
cyclotomic q-Schur algebras, Representation theory of algebraic
groups and quantum groups, Adv. Studies Pure Math. (2004): 261--320.

\end{thebibliography}
\end{document}